\documentclass[11pt]{article}
\usepackage{hyperref}
\pdfoutput=1
\usepackage[english]{babel}
\usepackage{amsmath,amsthm}
\usepackage{subeqnarray}
\usepackage{amsfonts}
\usepackage{graphicx}
\usepackage{float}
\usepackage[onehalfspacing]{setspace}
\usepackage[left=1.45in, top=1.25in, right=0.75in, bottom=1.15in]{geometry}

\theoremstyle{definition}

\theoremstyle{remark}

\numberwithin{equation}{section}

\newcommand{\set}[2]{\left\{{#1}\,:~{#2}\right\}}
\newcommand {\average}[1] {\mbox{$\left\{\!\!\left\{ #1 \right\}\!\!\right\}$}}
\newcommand {\jump}[1] {\mbox{$\left[\!\left[ #1 \right]\!\right]$}}

\begin{document}

\title{Local Improvements to Reduced-Order Approximations of Optimal Control Problems Governed by Diffusion-Convection-Reaction Equation}
\author{Tu\u{g}ba Akman \thanks{Electronic address: \texttt{takman@metu.edu.tr} This research was supported by the Middle East
Technical University Scientific Research Fund (Project: BAP-07-05-2013-003) }\\
\vspace{6pt} Department of Mathematics and Institute of Applied Mathematics \\
Middle East Technical University, 06800 Ankara, Turkey}

\maketitle

\begin{abstract}
  We consider the optimal control problem governed by diffusion convection reaction equation without control constraints. The proper orthogonal decomposition(POD) method is used to reduce the dimension of the problem. The POD method may be lack of accuracy if the POD basis depending on a set of parameters is used to approximate the problem depending on a different set of parameters. We are interested in the perturbation of diffusion term. To increase the accuracy and robustness of the basis, we compute three bases additional to the baseline POD. The first two of them use the sensitivity information to extrapolate and expand the POD basis. The other one is based on the subspace angle interpolation method. We compare these different bases in terms of accuracy and complexity and investigate the advantages and main drawbacks of them.
\end{abstract}


\section{Introduction}\label{introduction}

Optimal control problems(OCPs) governed by partial differential equations(PDEs) appear in several applications such as fluid dynamics, environmental modelling and engineering. To achieve an accurate numerical solution, a problem with a large dimension has to be solved. Even, it is resolved several times if an iterative optimization method is preferred. This increases the computational complexity. We use discontinuous Galerkin(DG) method with upwinding the convection term so that the continuity condition on the test and the trial functions are flexed, particularly along the edges of the triangles. However, the degrees of freedom is much higher than for the continuous Galerkin method. Thus, in order to solve the problem fast, we employ proper orthogonal decomposition(POD) method to derive the low-dimensional model of the problem.

POD is an efficient tool to produce an optimal basis from the numerical solution of the problem at discrete time steps. The POD basis, as opposed to the finite element method, contains information about the dominant characteristics of the problem. The accuracy of the reduced solution depends on how much information about the full-order solution is contained in this basis. More accurate solutions can be found by increasing the number of POD basis functions. However, this increases the computational cost. Therefore, the basis dimension is decided by balancing the accuracy and the number of the truncated POD basis. Specifically, the most energetic POD modes are chosen by measuring the ratio between the eigenvalues of the retained POD basis and the sum of the whole eigenvalues. 

In the literature, there are several studies concerning the reduced-order modelling. An overview of model reduction techniques for fluid dynamics problems can be found in the following report \cite{TLassila_AManzoni_AQuarteroni_2013a}. However, POD reduced-order modelling of the OCPs governed by diffusion-convection-reaction equation is rare, especially in 2D. For example, a finite horizon OCP for one-dimensional advection-diffusion equation is solved by choosing the time steps adaptively in \cite{AAlla_MFalcone_2013a}. In \cite{KKunisch_SVolkwein_1999a}, the reduced-order solution for the OCP governed by Burgers' equation is investigated. Optimality system-POD method is proposed in order to eliminate the effect of the reference control, which might be different from the optimal control, to the state in \cite{KKunisch_SVolkwein_2008a}. The choice of the time steps from which POD basis is computed can be decided in such a way that the error between the POD solution and its corresponding trajectory is minimized \cite{KKunisch_SVolkwein_2010a}. On the other hand, POD is applied to reacting flows with an adaptive strategy to capture the local dynamics of the chemical process in \cite{MSinger_HWGreen_2009a} in order to reduce the computational cost. In \cite{MLi_PDChristofidess_2008a}, optimal control of diffusion-convection-reaction process is concerned by applying two different spatial discretization techniques to compute the POD basis. In the presence of control or state constraints, a posteriori error estimation is employed to measure the difference between the suboptimal control computed using the POD basis and the optimal control and then the number of POD basis is decided \cite{FTroeltzsch_SVolkwein_2009a,AStudinger_SVolkwein_2013a,GMartin_SVolkwein_2014a}.

The POD basis which is generated via the snapshot ensemble depending on a set of parameters may not give accurate results when it is used to approximate the problem depending on a different set of parameters. For an accurate reduced solution of the perturbed problem, the POD basis must be regenerated for each parameter, which is costly. In this paper, motivated by the studies \cite{AHay_JTBorggaard_PDominique_2009a,AHay_JTBorggaard_IAkhtar_2010a,AHay_AImran_JTBorggaard_2012a}, we include the sensitivity information to enrich the low-dimensional subspace by extrapolating and expanding the POD basis in case of perturbation of the diffusion term. This idea requires the solution of the sensitivity equations which can be obtained by applying continuous sensitivity equations(CSEs) method or finite difference(FD) approximation. Because CSEs are always linear, the former method is especially preferable for nonlinear problems. The latter one requires the computation of the full problem at least one more time, so it is expensive for nonlinear case. Apart from the sensitivity analysis method, we use the subspace angle interpolation method(SAIM), which is applied to fluid-structure interaction problems in \cite{FVetrano_LCGarrec_2011a}, to increase the robustness. We construct the snapshot matrix using the state solution, the adjoint solution or a combination of them \cite{MHinze_SVolkwein_2008a}. We compare these different bases in terms of accuracy and complexity by measuring the error and the computational time and investigate the advantages and main drawbacks of them.

The rest of the paper is organized as follows: In Sect.\ref{ocp}, we introduce the OCP governed by diffusion-convection-reaction equation without control constraints. In Sect.\ref{prob_disct}, the symmetric interior penalty Galerkin method is explained and the fully-discrete optimality system is derived. The POD method and the reduced-order optimality system follow in Sect.\ref{pod}. In Sect.\ref{pod_basis}, the derivation of different POD bases using the sensitivity analysis and SAIM are explained. We present the numerical results in Sect.\ref{numerical}. Then, the paper ends with the conclusion.


\section{The Optimal Control Problem}\label{ocp}
In this paper, we are interested in the following distributed optimal control problem governed by the unsteady diffusion convection reaction equation without control constraints
\begin{align}\label{s1}
 \underset{u \in L^2(0,T; L^2(\Omega))} {\hbox{ minimize }} \; J(y,u):= \int_{0}^{T} \big( \frac{1}{2} \left\|y-y_{d}\right\|^{2}_{L^{2}(\Omega)} \ +
                    \frac{\alpha}{2}  \left\|u\right\|^{2}_{L^{2}(\Omega)} \big) \; dt,
\end{align}
subject to
\vspace{-4mm}
\begin{subequations}
\begin{align} 
\partial_t y-\epsilon \Delta y+ \boldsymbol{\beta} \cdot\nabla y+r y& =  f + u,  \;\qquad x\in \Omega, \;\;\quad t \in (0,T], \label{s2_a}  \\
  y(x,t)&=0,    \;\;\quad \quad \qquad x \in \partial \Omega, \quad t \in (0,T], \label{s2_b}  \\
 y(x,0)&=y_{0}(x), \quad \qquad  x\in \Omega, \label{s2_c}
\end{align}\label{s2}
\end{subequations}
where $\Omega$ is a bounded open, convex domain in  $\mathbb{R}^2$ with a Lipschitz boundary $\partial \Omega$ and $I = (0, T]$ is the time interval. $f, y_d \in L^{2}(0,T;L^{2}(\Omega)),\; y_0(x) \in H_{0}^1(\Omega), \boldsymbol{\beta} \in (W^{1,\infty}(\Omega))^2$ are given functions and $r, \epsilon, \alpha$ are given positive scalars. The velocity field $\boldsymbol{\beta}$ does not depend on time and satisfies the incompressibility condition, i.e. $\nabla \cdot \boldsymbol{\beta} =0$. 

In order to write the variational formulation of the problem, we define the bilinear forms 
\begin{align}
a(y,v)=\int_{\Omega} (\epsilon \nabla y \cdot \nabla v + \boldsymbol{\beta} \cdot \nabla y v + r y v)\; dx, \quad 
(u,v)=\int_{\Omega} uv \; dx,
\end{align}
the state and the test space as $Y=V=H^1_0(\Omega), \forall t \in (0,T]$. It is well known that the pair $(y,u) \in H^{1}(0,T;L^{2}(\Omega)) \cap L^{2}(0,T;H^1_0(\Omega)) \times L^2(0,T; L^{2}(\Omega))$ is the unique solution of the optimal control problem (\ref{s1}-\ref{s2}) if and only if there is an adjoint $p \in H^{1}(0,T;L^{2}(\Omega)) \cap L^{2}(0,T;H^1_0(\Omega))$ such that $(y,u,p)$ satisfy the following optimality system \cite{FTroeltzsch_2010a}
\begin{subequations}
\begin{eqnarray}
(\partial_t y,v)+a(y,v) & = & (f + u,v), \qquad \forall v \in V,  \; \; \; y(x,0)=y_0, \label{s4_a} \\
 -(\partial_t p,q)+a(q,p) & = & -(y-y_{d},q), \quad \forall q \in V, \; \; p(x,T)=0, \label{s4_b} \\
\alpha u &=& p. \label{s4_c}
\end{eqnarray}\label{s4}
\end{subequations}


\section{Problem Discretization}\label{prob_disct}
In this section, we briefly describe the symmetric interior penalty Galerkin(SIPG) method which we use for spatial discretization.
Let $\{ \mathcal{T}_h\}_h$ be a family of shape regular meshes such that $\overline{\Omega} = \cup_{K \in \mathcal{T}_h} \overline{K}$, $K_i \cap K_j = \emptyset$ for $K_i, K_j \in \mathcal{T}_h$, $i \not= j$. The diameter of an element $K$ and
the length of an edge $E$ are denoted by  $h_{K}$  and $h_E$, respectively. Further, the maximum value of element diameter
is denoted by $h=\max \limits_{K \in \mathcal{T}_h} h_{K}$.

We only consider  discontinuous piecewise finite element spaces to
define the discrete spaces of the state, adjoint, control and test functions
\begin{align*}\label{Yspace}
V_h = Y_h = U_h &= \set{y \in L^2(\Omega)}{ y\mid_{K}\in \mathbb{P}^p(K) \quad \forall K \in \mathcal{T}_h}.
\end{align*}
Here, $\mathbb{P}^p(K)$ denotes the set of all polynomials on $ K \in \mathcal{T}_h$ of degree $p$. We split the  set of all edges $\mathcal{E}_h$ into the set $\mathcal{E}^0_h$ of interior edges  and the set $\mathcal{E}^{\partial}_h$ of boundary edges so that $\mathcal{E}_h=\mathcal{E}^{\partial}_h\cup \mathcal{E}^{0}_h$. Let $\mathbf{n}$ denote the unit outward normal to $\partial \Omega$. We define the inflow boundary
\[
          \Gamma^- = \set{x \in \partial \Omega}{\boldsymbol{\beta}\cdot \mathbf{n}(x) < 0}
\]
and the outflow boundary $\Gamma^+ = \partial \Omega \setminus  \Gamma^-$. The boundary edges are decomposed into edges
$\mathcal{E}^{-}_h = \set{ E \in \mathcal{E}^{\partial}_h}{ E \subset  \Gamma^- }$ that correspond to inflow boundary and edges
$\mathcal{E}^{+}_h = \mathcal{E}^\partial_h \setminus \mathcal{E}^{-}_h$ that  correspond to outflow boundary. The inflow and outflow boundaries of an element $K \in \mathcal{T}_h$ are defined by
\begin{align*}
\partial K^-&=\set{x \in \partial K}{\boldsymbol{\beta} \cdot \mathbf{n}_{K}(x) <0}, \quad \partial K^{+} = \partial K \setminus \partial K^{-},
\end{align*}
where $\mathbf{n}_{K}$ is the unit normal vector on the boundary $\partial K$ of an element $K$.

Let the edge $E$ be a common edge for two elements $K$ and $K^e$. For a piecewise continuous scalar function $y$,
there are two traces of $y$ along $E$, denoted by $y|_E$ from inside $K$ and $y^e|_E$ from inside $K^e$. Then, the jump and average of $y$ across the edge $E$ are defined by:
\begin{align*}
\jump{y}&=y|_E\mathbf{n}_{K}+y^e|_E\mathbf{n}_{K^e}, \quad
\average{y}=\frac{1}{2}\big( y|_E+y^e|_E \big).
\end{align*}

Similarly, for a piecewise continuous vector field $\nabla y$, the jump and average across an edge $E$ are given by
\begin{align*}
\jump{\nabla y}&=\nabla y|_E \cdot \mathbf{n}_{K}+\nabla y^e|_E \cdot \mathbf{n}_{K^e}, \quad
\average{\nabla y}=\frac{1}{2}\big(\nabla y|_E+\nabla y^e|_E \big).
\end{align*}

For a boundary edge $E \in K \cap \Gamma$, we set $\average{\nabla y}=\nabla y$ and $\jump{y}=y\mathbf{n}$
where $\mathbf{n}$ is the outward normal unit vector on $\Gamma$.

We can now give DG discretizations of the state equation (\ref{s2}) in space  for fixed control $u$.  The DG method proposed here is based on the upwind discretization of the convection term and on the SIPG discretization  of the diffusion term \cite{DSchotzau_LZhu_2009a}. This leads to the following (bi-)linear forms applied to $y_{h} \in H^1(0,T;Y_h)$ for $\forall t \in (0,T]$
\vspace{-2mm}
\begin{align*}
(\partial_t y_h, v_h) + a_h(y_h,v_h) =(f_h + u_h,v_h)  \quad \forall v_h \in V_h, \quad  t \in (0,T],
\end{align*}
where
\vspace{-3mm}
\begin{align*} 
a_h(y,v)=& \sum \limits_{K \in \mathcal{T}_h} \int \limits_{K} \epsilon \nabla y \cdot  \nabla v \; dx           
         +  \sum \limits_{K \in \mathcal{T}_h} \int \limits_{K} \big( \boldsymbol{\beta} \cdot \nabla y v + r y v \big) \; dx \nonumber \\
         &- \sum \limits_{ E \in \mathcal{E}_h} \int \limits_E \big( \average{\epsilon \nabla y} \cdot \jump{v}
                    + \average{\epsilon \nabla v} \cdot \jump{y} - \frac{\sigma \epsilon}{h_E}  \jump{y} \cdot \jump{v} \big) \; ds \nonumber \\
         &+ \sum \limits_{K \in \mathcal{T}_h}\; \int \limits_{\partial K^{-} \backslash \Gamma^-} \boldsymbol{\beta} \cdot \mathbf{n} (y^e-y)v \; ds
           - \sum \limits_{K \in \mathcal{T}_h} \; \int \limits_{\partial K^{-} \cap \Gamma^{-}} \boldsymbol{\beta} \cdot \mathbf{n} y v  \; ds,
\end{align*}
with a constant  interior penalty parameter $\sigma >0$. We choose $\sigma$ to be sufficiently large, independent of the mesh size $h$ and the diffusion coefficient $\epsilon$ to ensure the stability of the DG discretization as described in \cite[Sec.~2.7.1]{BRiviere_2008a} with a lower bound depending  only on the polynomial degree.

Let $f_h, y_h^d$ and $y_h^0$ be approximations of the source function $f$, the desired state function $y_d$ and initial condition $y_0$, respectively. The functions $(y_{h},u_{h}) \in H^1(0,T;Y_h) \times L^2(0,T;U_{h})$ solve the semi-discrete OCP if and only if $(y_h,u_{h},p_h) \in H^1(0,T;Y_h) \times L^2(0,T;U_{h}) \times H^1(0,T;Y_h)$ is a unique solution of the following optimality system:
\begin{subequations}\label{s7}
\begin{align}
(\partial_t y_h, v_h) + a_{h}(y_{h},v_{h})+b(u_{h},v_{h})&=(f_h,v_{h}),  & \forall v_{h} \in V_{h}, \quad y_h(x,0)=y_h^0,\\ 
-(\partial_t p_h, q_{h}) + a_{h}(q_{h},p_{h})&=-(y_{h}-y_h^{d},q_{h}), & \forall q_{h} \in V_{h}, \quad p_h(x,T)=0,\\ 
\alpha u_{h} &=p_{h}.
\end{align}
\end{subequations}
For details, we refer the reader to the study \cite{TAkman_HYucel_BKarasozen2013a} where a priori error estimates for SIPG discretization combined with backward Euler is provided and the quadratic convergence rate in space is achieved. For time-discretization, we use Crank-Nicolson method, which is known to be stable and second order convergent. 

Let $0 = t_0 < t_1 < \cdots < t_{N} = T $ be a subdivison of $I$ with time intervals $I_m = (t_{m-1}, t_{m}]$ and time steps $k = t_{m} - t_{m-1}$ for $m=1, \ldots, N$. Then, the fully discrete optimality system for (\ref{s7}) is written as follows:
\begin{subequations}\label{OD_ocp_theta}
\begin{align}
 (y_{h,m+1}-y_{h,m},v_{h}) + \frac{k}{2}a_{h}(y_{h,m}+ y_{h,m+1},v_{h}) &= 
 \frac{k}{2} ( (f_{h,m} + u_{h,m},v_{h}) + (f_{h,m+1} + u_{h,m+1},v_{h})), \\
 m&=0, \cdots, N-1,  \qquad
 y_{h,0}(x,0) =y_{0}, \nonumber \\
 (p_{h,m} - p_{h,m+1},q_{h}) + \frac{k}{2} a_{h}(q_{h},p_{h,m}+p_{h,m+1}) &=  
 -\frac{k}{2}\left((y_{h,m} - y^{d}_{h,m},q_{h}) + (y_{h,m+1}-y^{d}_{h,m+1},q_{h})\right), \\
 m &= N-1, \cdots, 0, \qquad
 p_{h,N} = 0, \nonumber \\
 \alpha u_{h,m} &=p_{h,m}, \quad m=0,1, \ldots, N. 
\end{align}
\end{subequations}
\section{Proper Orthogonal Decomposition}\label{pod}
In this section, we briefly explain the POD method following \cite{SVolkwein_2008a} and its generalization to Galerkin type spatial discretization \cite[Sec. 3]{MFahl_2000a}. Let the DG approximation to the function $w(x,t_i)$ at $t=t_i$ be expressed as
\begin{equation*}
w_i^{DG} = \sum\limits_{j=1}^{M} \textbf{w}^{DG}_j(t_i) \varphi_j(x),
\end{equation*}
where $\varphi_j(x)$ denotes the discontinuous finite element functions with $x \in \Omega$. Let the matrix $W$ be a real-valued $M\times N$ matrix containing the DG coefficients $\textbf{w}^{DG}_j(t_i)$ on its $i$-th column. We assume that the snapshots are linearly independent so that $W$ has full-column rank. On the other hand, the POD approximation is written as
\begin{equation*}
w_i^{POD} = \sum\limits_{j=1}^{l} \textbf{w}^{POD}_j(t_i) \psi_j(x),
\end{equation*}
where $\psi_j(x)$ denotes the POD basis functions. The aim is to find the function $ \psi_j(x)$ representing the snapshot ensemble $w_i^{DG}$ as well as possible, which is equivalent to the following minimization problem \cite{KKunisch_SVolkwein_1999a}
\begin{equation*}
\min \sum \limits_{i=1}^{M} \left\|  w_i^{DG} -  \sum \limits_{j=1}^{l}(w_i^{DG}, \psi_j)_{\mathcal{M}}  \psi_j \right\|^2_{\mathcal{M}} \ \hbox{s.t. } \ \| \psi\|_{\mathcal{M}} = 1,
\end{equation*}
where $(\phi, \psi)_{\mathcal{M}} =  \boldsymbol{\phi}^{T} \mathcal{M} \boldsymbol{\psi}$ with the finite element mass matrix $\mathcal{M}$. 

The POD basis is computed using the singular value decomposition(SVD), because it is more stable than the eigenvalue decomposition, i.e. the singular values decay to machine precision, whereas the eigenvalues stagnate above \cite{AStudinger_SVolkwein_2013a}. The Cholesky decomposition of the symmetric positive definite mass matrix $\mathcal{M}~=~\mathcal{M}^{1/2}~(\mathcal{M}^{1/2})^T$ is required. Then, we obtain the data matrix $\tilde{W}~ =~ M^{1/2} W$ and whose SVD is written as $\tilde{W}~ =~ U \Sigma V^T$. Firstly, coefficients  of the POD basis, namely $\Psi$, are computed by solving the linear system 
\begin{equation}\label{pod_Psi1}
(M^{1/2})^T \Psi = U^l,
\end{equation}
with the first $l$ columns of $U$. Then, $l$-many POD basis functions are expressed as a linear combination of the finite element basis functions $\varphi_{i}(x)$ as follows
\begin{equation}\label{Psi1}
\psi_j(x) = \sum \limits_{i=1}^{M} \Psi_{ij} \varphi_{i}(x), \quad j=1,\ldots, l.
\end{equation}

The dimension of the POD basis must be several order lower than the dimension of the full-order space from which we obtain the snapshots of the problem. The number of POD basis functions $l$ is decided according to the ratio between the modelled and the total energy,
\begin{equation}\label{energy}
\mathcal{E}(l) = \sum \limits_{i=1}^{l} \sigma_i^{2} / \sum \limits_{i=1}^{d} \sigma_i^{2},
\end{equation}
where $\sigma_i$'s denote the singular values of the data matrix $\tilde{W}$ and $d=\mbox{rank}(\tilde{W})$.

After deriving the POD basis, low-dimensional optimality system corresponding to (\ref{OD_ocp_theta}) is derived by applying Galerkin projection. In addition, the initial condition $y_0$, desired state $y_d$ and the source function $f$ are also projected onto the low-dimensional space. Then, the reduced-order POD Galerkin modelling of the optimality system is written as follows:
\begin{align*}\label{OD_POD}
 (y_{m+1}^{POD}-y_{m}^{POD},\psi_j) + \frac{k}{2}a_{h}(y_{m}^{POD}+ y_{m+1}^{POD},\psi_j) &= 
 \frac{k}{2} ( (f_{m}^{POD} + u_{m}^{POD},\psi_j) + (f_{m+1}^{POD} + u_{m+1}^{POD},\psi_j)), \nonumber \\
 m&=0, \cdots, N-1,  \qquad
 y_{0}^{POD}(x,0) =y_{0}, \nonumber \\
 (p_{m}^{POD} - p_{m+1}^{POD},\psi_{j}) + \frac{k}{2} a_{h}(\psi_{j},p_{m}^{POD}+p_{m+1}^{POD}) &=  
 -\frac{k}{2}\left((y_{m}^{POD} - y^{d,POD}_{m},\psi_{j}) + (y_{m+1}^{POD}-y^{d,POD}_{m+1},\psi_{j})\right),\nonumber \\
 m &= N-1, \cdots, 0, \qquad
 p_{N}^{POD} = 0, \nonumber \\
 \alpha u_{m}^{POD} &=p_{m}^{POD}, \quad m=0,1, \ldots, N. \nonumber
\end{align*}

\section{Computation of Different POD Bases}\label{pod_basis}

In  general, the POD basis generated via the snapshots depending on a parameter $\mu_0$ cannot capture the dynamics of the perturbed problem depending on $\mu = \mu_0 + \Delta \mu$. In order to eliminate this drawback of the method, motivated by the studies \cite{AHay_JTBorggaard_PDominique_2009a,AHay_JTBorggaard_IAkhtar_2010a,AHay_AImran_JTBorggaard_2012a}, POD sensitivities can be used to enrich the low-dimensional space for a wider range of parameters. We generate two new bases, i.e. extrapolated POD(ExtPOD) and expanded POD(ExpPOD). In addition to these, the subspace angle interpolation method(SAIM) \cite{ABjork_GHGolub_1973a} can be applied, which does not require the sensitivity information. Now, we proceed with the details of the sensitivity analysis.
 
\subsection{Proper Orthogonal Decomposition Sensitivities}\label{sens}
The sensitivity of a term is defined as the derivative of that term with respect to a quantity of interest $\mu$. We assume that the state, the adjoint and the control are functions depending on space, time and the quantity of interest $\mu$,  
\[
y=y(x,t,\mu),\quad p=p(x,t,\mu), \quad u=u(x,t,\mu).
\] 
We define the sensitivities as
\begin{equation*}
s_y = \frac{ \partial y}{ \partial \mu}, \quad 
s_p = \frac{ \partial p}{ \partial \mu}, \quad 
s_u = \frac{ \partial u}{ \partial \mu}.
\end{equation*}

In this study, we are interested in the sensitivity with respect to the diffusion term $\mu = \epsilon$. The CSEs are obtained by differentiating the continuous state (\ref{s2}) and the adjoint equation associated to (\ref{s4_b}) and the optimality condition (\ref{s4_c}) with respect to $\mu$. The subscript $\mu$ denotes the derivative with respect to $\mu$. The corresponding optimality system with $s_y, s_p$ and $s_u$ is written as follows,
\begin{align}\label{sens_s1}
(\partial_t s_y,v)+a(s_y,v) + (\nabla y, \nabla v)&=(f_{\mu} + s_u,v),\qquad \forall v \in V,  \; \; \;  s_y(x,0)=(y_0)_{\mu},  \nonumber\\
-(\partial_t s_p,\psi)+a(\psi, s_p) + (\nabla \psi, \nabla p)&=-(s_y-y_{\mu}^d,\psi), \quad \forall \psi \in V,  \; \; \;   s_p(x,T)=0,  \\
\alpha s_u &= s_p.\qquad \qquad \; \;  \nonumber
\end{align}


We note that the homogeneous Dirichlet boundary conditions are imposed to (\ref{sens_s1}) after differentiating (\ref{s2_b}) in the same way. The optimality system (\ref{sens_s1}) is discretized using the same numerical method, i.e. SIPG in space and Crank-Nicolson time, as for (\ref{s4}). 

The sensitivity equations are always linear, so CSE method would be especially promising for nonlinear problems. On the other hand, FD approximation can also be used to find the sensitivities. It requires the evaluation of the OCP depending on the perturbed parameters. In particular, the sensitivity of the state can be computed via the centred difference
\begin{equation}\label{FD}
s_{y}(\mu_0) \approx \frac{y(\mu_0 + \Delta \mu) - y(\mu_0 - \Delta \mu)}{2 \Delta \mu}.
\end{equation}

The increment $ \Delta \mu$ is chosen sufficiently small for an accurate FD approximation and it is chosen sufficiently large for the difference between two nearby POD vectors to be larger than the discretization error by one order of magnitude \cite{AHay_JTBorggaard_PDominique_2009a}.

After finding the sensitivities of the state and the adjoint, the POD sensitivities are obtained. To do so, we treat each POD mode as a function of both space and the parameter, i.e. $\psi = \psi(x, \mu)$. Then, we differentiate the relation (\ref{pod_Psi1}) with respect to $\mu$ and solve the resulting equation for $\Psi_{\mu}$. We proceed with the relation (\ref{Psi1}) to derive the POD basis sensitivities $\psi_{\mu}$
\begin{equation*}
{(\psi_j)}_{\mu} = \sum \limits_{i=1}^{m} {(\Psi_{ij})}_{\mu} \varphi_{i}(x), \quad j=1,\ldots, l.
\end{equation*}

The computation of $U^l_{\mu}$, which appears after differentiating (\ref{pod_Psi1}) with respect to $\mu$, is realised through the relation
\begin{equation*}
U^l_{\mu}
= (\tilde{W} V^{l} \Sigma^{\dagger})_{\mu}
= \tilde{W}_{\mu} V^{l} \Sigma^{\dagger} + \tilde{W} V^{l}_{\mu}  \Sigma^{\dagger} + \tilde{W} V^{l} \Sigma^{\dagger}_{\mu}.
\end{equation*}

The term $\tilde{W}_{\mu}$ denotes the sensitivity of the snapshot matrix, which is obtained by solving (\ref{sens_s1}) using CSE method or FD approximation. For the computation of $V_{\mu}^{l}$ and $\Sigma^{\dagger}_{\mu}$, we consider the eigenvalue problem $B V^l = V^l \lambda^l$ with the $l$th column of $V$. 
After differentiation, one obtains the relation
\begin{equation}\label{VBV}
(V^l)^T (B_{\mu} - \lambda^l_{\mu} I)V^l = 0.
\end{equation}
The orthonormal matrix $V^l$ has already been computed via SVD. Then, the eigenvalue sensitivities are given by
\begin{equation*}\label{eig_sens}
\lambda_{\mu}^{l} = (V^l)^T B_{\mu} V^l.
\end{equation*}
Each term of $ \Sigma^{\dagger}_{\mu}$ is computed due to the relation between the singular values $\sigma_{\mu}$ and the eigenvalues $\lambda_{\mu}$, i.e. $\sigma^2_{\mu} = \lambda_{\mu}$. 

The equation (\ref{VBV}) is solved in the least-squares sense and we denote one particular solution by $s^l$. The general solution to (\ref{VBV}) is expressed as $s^{l} + \gamma V^{l}$ for $\gamma \in \mathbb{R}$ with a simple $\lambda^{l}$. In addition, we differentiate the normalization condition $V^l (V^l)^T=~1$ leading to $V^l_{\mu} (V^l)^T=0$. Then, the sensitivity of $V^l$ and $\gamma$ are determined by
\begin{equation*}
V^l_{\mu} = s^l - ((s^l)^T V^l)V^l, \qquad \gamma = -(s^l)^T V^l.
\end{equation*} 
For details, we refer the reader to \cite[Sec. 3.2]{AHay_JTBorggaard_PDominique_2009a}. 

In ExtPOD, the POD basis depending on $\mu$ is written using the first-order Taylors expansion as follows
\begin{equation*}
\psi(x, \mu) = \psi(x, \mu_0) + \Delta \mu\frac{\partial \psi}{\partial \mu}(x, \mu_0) + \mathcal{O}(\Delta \mu^2).
\end{equation*}
The reduced-order solution is expressed as
\begin{equation*}
w_i^{POD} = \sum \limits_{j=1}^{l} \textbf{w}^{POD}_j(t_i) (\psi_j(x, \mu_0) + \Delta \mu (\psi(x, \mu_0)_{\mu}).
\end{equation*}

In ExpPOD, the POD basis sensitivities are also added to the original POD basis as follows
$$[\psi_1, \ldots, \psi_l, ({\psi_{1}})_{\mu}, \ldots, ({\psi_l})_{\mu}]$$ 
and the reduced-order solution is written as
\begin{equation*}
w_i^{POD} = \sum \limits_{j=1}^{l} \textbf{w}^{POD}_j(t_i) \psi_j(x, \mu_0) + \sum \limits_{j=l+1}^{2l} \textbf{w}^{POD}_j(t_i) (\psi_{j-l}(x, \mu_0))_{\mu},
\end{equation*}
where the dimension of the reduced basis is doubled.

\subsection{The Subspace Angle Interpolation Method}\label{saim}
In this section, we explain the method following \cite[Thm.1]{ABjork_GHGolub_1973a} and \cite{FVetrano_LCGarrec_2011a}. Let $\Psi^{1}$ and $\Psi^{2}$ be the coefficients of two POD bases spanning two subspaces associated to two different parameters $\mu_1$ and $\mu_2$. We derive SVD of their product, which is given by $(\Psi^{1})^{T}\Psi^{2} = \tilde{U} \tilde \Sigma \tilde{V}^{T}$. For a decreasing sequence of the singular values of $\tilde \Sigma$, i.e. $\tilde \sigma_1 \geq \tilde \sigma_2 \geq \ldots \geq \tilde \sigma_q \geq 0$, the principal angles and the principal vectors associated to those subspaces are given by
\begin{equation*}
\cos(\theta) = \tilde \Sigma, \qquad U = \Psi^{1}\tilde{U},\quad V = \Psi^{2}\tilde{V}^{T}.
\end{equation*} 

Then, the coefficients of the new basis are constructed as in \cite[Sect.3.2]{FVetrano_LCGarrec_2011a} for any parameter  $\mu_1 < \mu_N < \mu_2$ through the linearly interpolated principle angle 
\begin{equation*}
\theta_j(\mu_1, \mu_N) = \left( \frac{\mu_N - \mu_1}{\mu_2 - \mu_1} \right) \theta_j(\mu_1, \mu_2), \quad j=1,\cdots,l.
\end{equation*} 

The coefficients of the new basis are written as
\begin{equation*}
\Psi_{:,j} = u_{:,j} \cos(\theta_j(\mu_1, \mu_N)) + \frac{v_{:,j} - (u_{:,j}^T v_{:,j})u_{:,j}}{|| v_{:,j} - (u_{:,j}^T v_{:,j})u_{:,j} ||^2} \sin(\theta_j(\mu_1, \mu_N)), \quad j=1,\cdots,l,
\end{equation*}
where $u_{:,j}, v_{:,j}$ are the columns of $U$ and $V$. We note that this method only requires the coefficients of two precomputed bases and it is applied if only one parameter is perturbed in the system. 

\section{Numerical Results}\label{numerical}
In this section, we present some numerical results to investigate the performances of the methods. The CPU times are obtained on a 3.17 GHz desktop PC. The full problem is solved with linear finite elements for $\Delta t=1/60$ and $\Delta x=1/40$ leading to 9600 degrees of freedom. We used the Newton-conjugate gradient method with Armijo line-search in the optimization step for fast convergence. Three different snapshot sets for $W$ are used to generate the POD basis, namely the state $Y$, the adjoint $P$ and the combination of them $Y \cup P$, as in \cite{MHinze_SVolkwein_2008a}. The sensitivities derived from CSE method are calculated at the same time steps with FD approximation and we use them in the bases generation step. The error between the benchmarked and the reduced solution is measured with respect to $L^2(0,T; L^2(\Omega))$ norm. 

We note that the nominal value for the diffusion term is $\epsilon=10^{-2}$. We generate POD basis once using the snapshots associated to this nominal/baseline value and denote the corresponding results by POD in the figures. We choose the parameter range for $\mu = \epsilon$ as $1/\epsilon = 80:5:120$. For SAIM, two snapshot sets associated to two different parameters are required. To do so, we fix the parameters as $\mu_1 = 1/125$ and $\mu_2=1/75$.

We consider the optimal control problem with
 \[
Q=(0,1] \times \Omega, \; \Omega=(0,1)^{2}, \; \epsilon=10^{-2}, \; \boldsymbol\beta=(y-1/2,-x+1/2)^T, \; r=1, \; \alpha=1.
\]
We take the source function $f$, the desired state $y_{d}$ and the initial condition $y_0$ as
\begin{align*}
f(x,t) &= y_d(x,t) = 1,          \\
y_0(x,t) &= 0.                
\end{align*}

\begin{figure}[H]
   $
   \begin{array}{ccc}
   \includegraphics[width=0.32\textwidth]{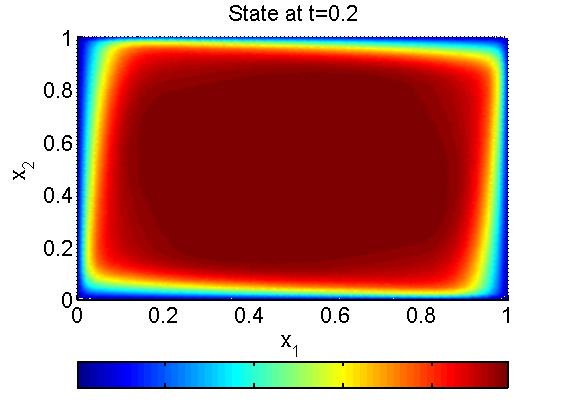}   
   \includegraphics[width=0.32\textwidth]{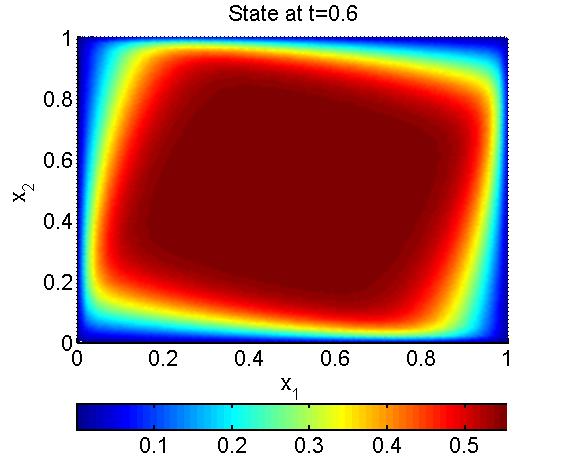}
   \includegraphics[width=0.32\textwidth]{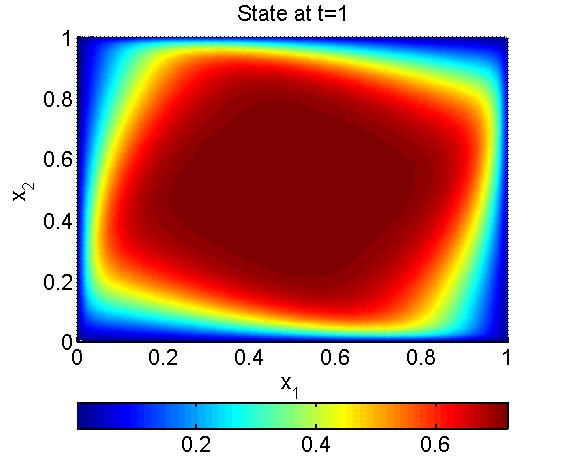}\\   
   \includegraphics[width=0.32\textwidth]{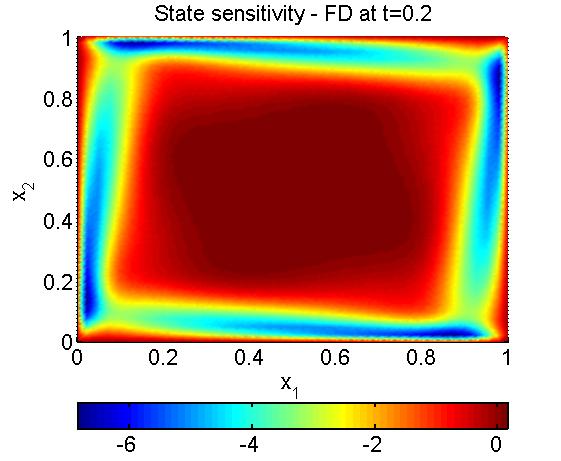}
   \includegraphics[width=0.32\textwidth]{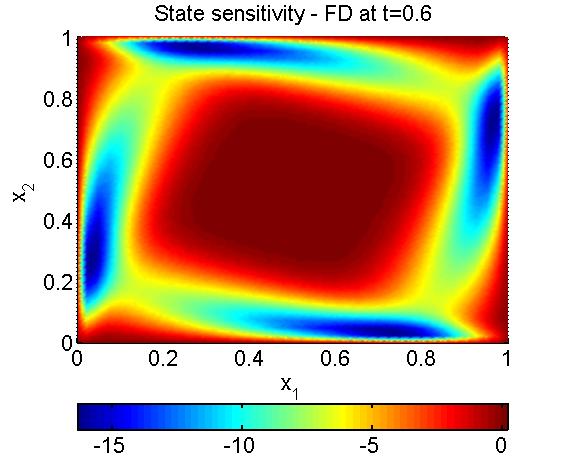}                          
   \includegraphics[width=0.32\textwidth]{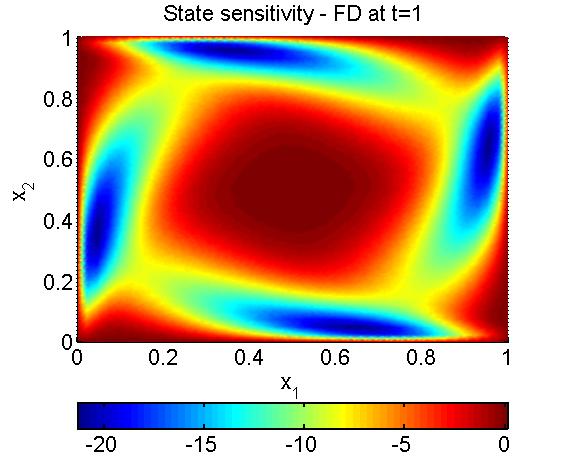}\\
   \end{array}
   $
   \caption{State(\textit{first row}) and state sensitivities(\textit{second row}) at $t=0.2, 0.6, 1$, respectively}
   \label{EX2_F:1}
\end{figure}

\begin{figure}[H]
   $
   \begin{array}{ccc}     
   \includegraphics[width=0.32\textwidth]{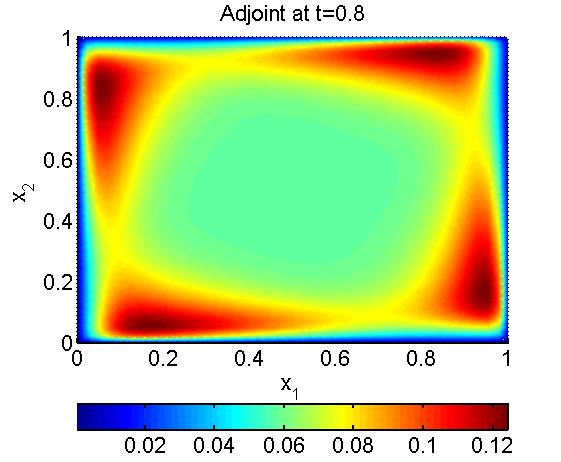}
   \includegraphics[width=0.32\textwidth]{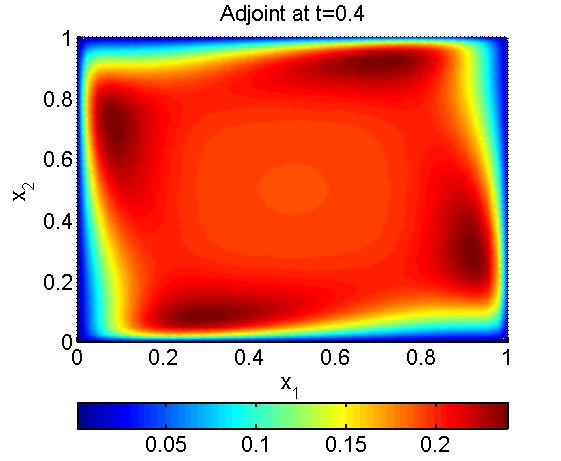}
   \includegraphics[width=0.32\textwidth]{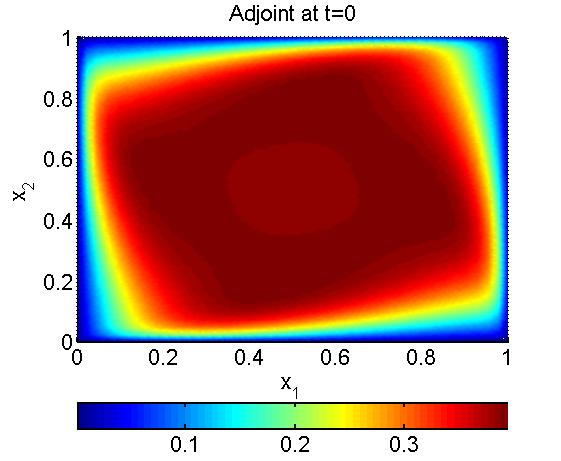}\\
   \includegraphics[width=0.32\textwidth]{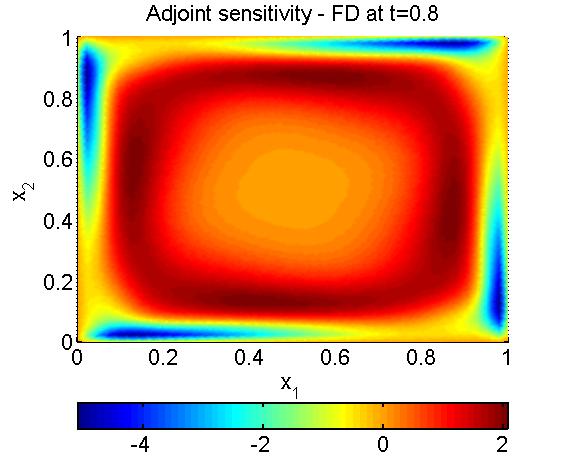}                
   \includegraphics[width=0.32\textwidth]{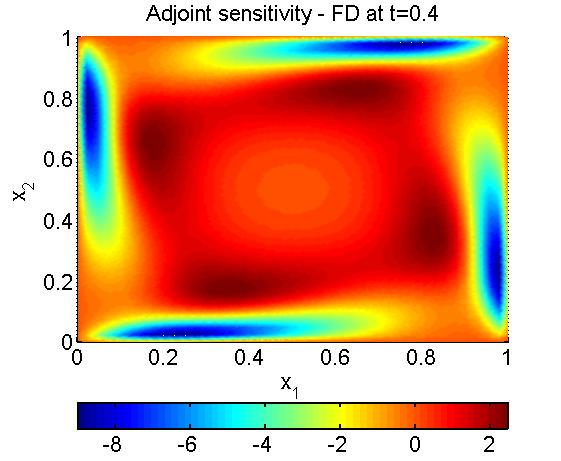}    
   \includegraphics[width=0.32\textwidth]{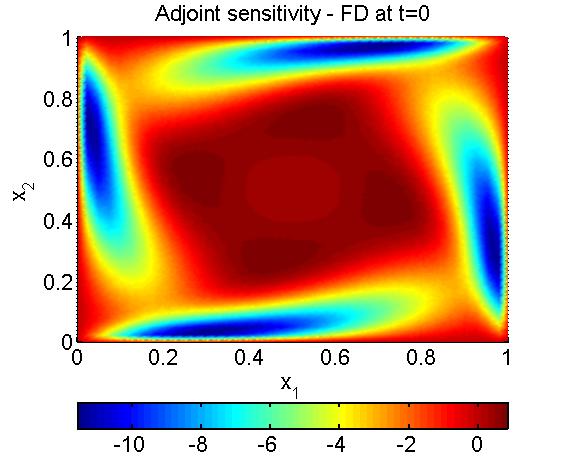}
   \end{array}
   $
   \caption{Adjoint(\textit{first row}) and adjoint sensitivities(\textit{second row}) at $t=0.8, 0.4, 0$, respectively}
   \label{EX2_F:2}
\end{figure}

The exact solution of this problem is not known and the convection field is not a constant vector. In Figure~\ref{EX2_F:1}, we present the numerical solutions of the state on the first row and its sensitivities on the second row at the time instances $t=0.2, 0.6, 1$ from left to right. In Figure~\ref{EX2_F:2}, the numerical solution of the adjoint on the first row and its sensitivities on the second row are given at the time instances $t=0.8, 0.4, 0$ from left to right. Due to the convection field $\boldsymbol\beta$, the state rotates clockwise as $t \longrightarrow T=1$, while the adjoint follows the counter-clockwise direction as $t \longrightarrow t_0=0$. Since the convection field is fixed, the associated sensitivities rotate in the same direction, too. The state and the adjoint are highly sensitive where their solutions change mostly. 

\begin{figure}[H]
   \centering
   $
   \begin{array}{cc}
   \includegraphics[width=0.48\textwidth]{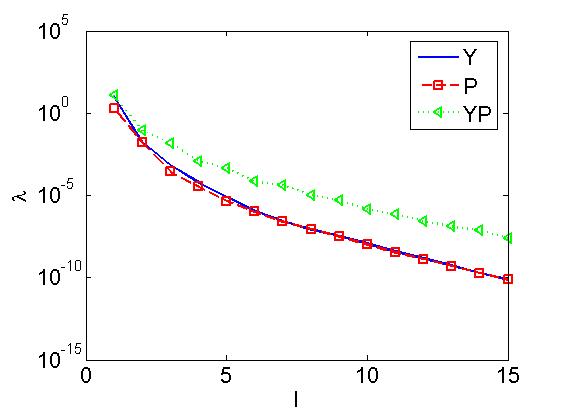}&
   \includegraphics[width=0.48\textwidth]{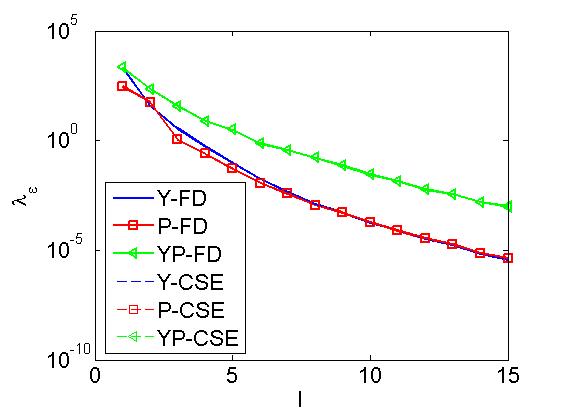}
   \end{array}
   $
   \caption{Eigenvalues(\textit{left}) and their sensitivities(\textit{right})}
   \label{EX2_F:3}
\end{figure}

In Figure~\ref{EX2_F:3}, we present the decay of the eigenvalues and their sensitivities. They decrease rapidly showing that POD can be successfully applied. The sensitivities are computed by CSE method and FD approximation for comparison purposes. We observe that both approaches give almost the same results. In addition, the eigenvalues decay following the same pattern as the sensitivities do which means that the ordering will remain in case of parameter perturbations \cite{AHay_JTBorggaard_PDominique_2009a}. 

We present the computational time for the full-problem, CSE method and FD approximation in Table~\ref{EX2_T:1}. In total, the full-problem is solved in 30 seconds. The CSE method takes almost 2 minutes. The sensitivities are computed using FD approximation in 55 seconds, which is faster than CSE method. 

\begin{table}[H]
\caption{CPU times for the full problem}
\label{EX2_T:1}
\begin{tabular}{lll}
\hline\noalign{\smallskip}
Computing the FE mesh and matrices   & $\approx$& 5 s \\
FE element solution                  & $\approx$& 25 s\\
Sensitivity solution using CSE method& $\approx$& 2 m\\
Sensitivity solution using FD approximation&$\approx$ & 2 $\times$ 25 s \\
\noalign{\smallskip}\hline
\end{tabular}
\end{table}

The total energy $\mathcal{E}(l)$ in the formula (\ref{energy}) is fixed up to $100(1-\gamma)\%$ by keeping the most energetic POD modes. In this study, we choose 9 POD basis functions setting $\gamma=10^{-2}$. In Table~\ref{EX2_T:2}, we compare the computational cost of the reduced problem in terms of different bases and the snapshot ensemble. For each case, the reduced problem is solved less than 5 seconds, which is faster than for the full-problem. Using the snapshot set $P$, the problem is solved faster than with the snapshot set $Y$. It is because a better approximation to the control is achieved and led to fast convergence in the optimization step. On the other hand, the size of the set $Y \cup P$ is twice as large as $Y$ or $P$. Therefore, it takes longer to compute the POD basis and the reduced solution. In terms of POD sensitivities, ExpPOD is slower than ExtPOD; because, its dimension is doubled. The SAIM is more costly than ExtPOD due to the additional computations to generate the new basis. We note that the speed of POD gains importance when we have to solve the full-problem several times in case of parameter perturbations.

\begin{table}[H]
\caption{CPU times in seconds for the reduced problem}
\label{EX2_T:2}
\begin{tabular}{lllll}
\hline\noalign{\smallskip}
 & & $Y$ & $P$ & $ Y \cup P $  \\
\noalign{\smallskip}\hline\noalign{\smallskip}
BPOD   &$\approx$& 1.48& 1.54  & 1.75\\
ExtPOD &$\approx$& 3.83& 3.13  & 3.88\\
ExpPOD &$\approx$& 4.68& 3.94  & 4.73\\
SAIM   &$\approx$& 4.07& 3.39  & 4.42\\
\noalign{\smallskip}\hline
\end{tabular}
\end{table}

\begin{figure}[H]
   \centering
      $
      \begin{array}{cc}  
   \includegraphics[width=0.48\textwidth]{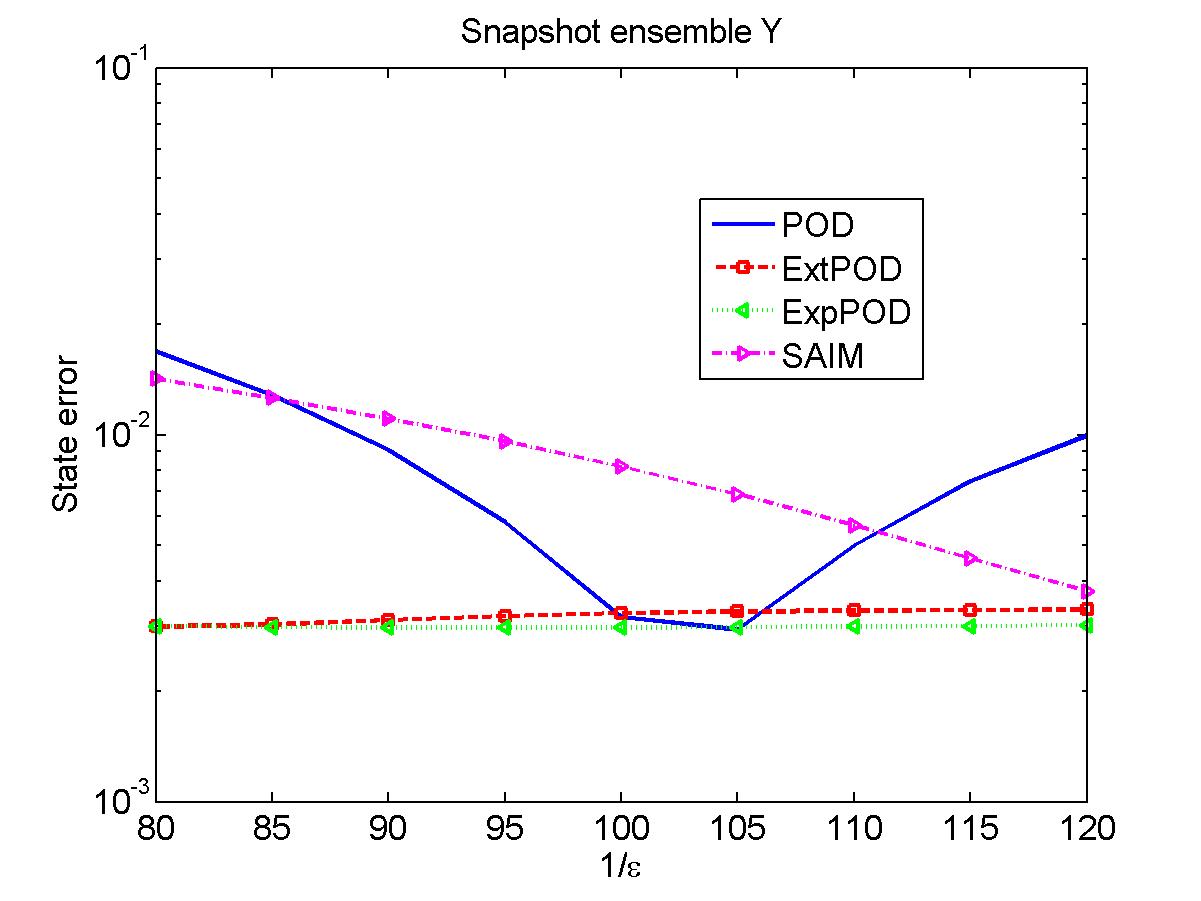}
   \includegraphics[width=0.48\textwidth]{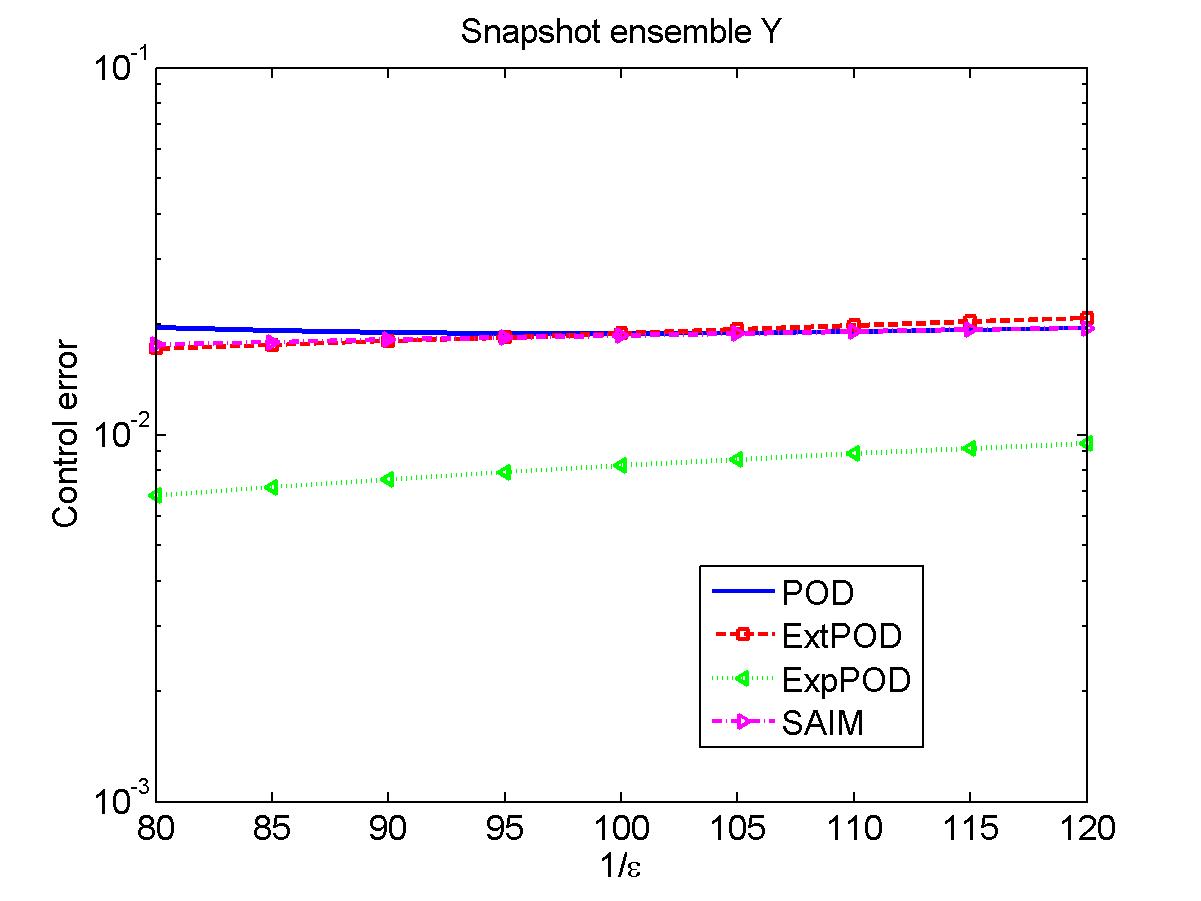} \\
   \includegraphics[width=0.48\textwidth]{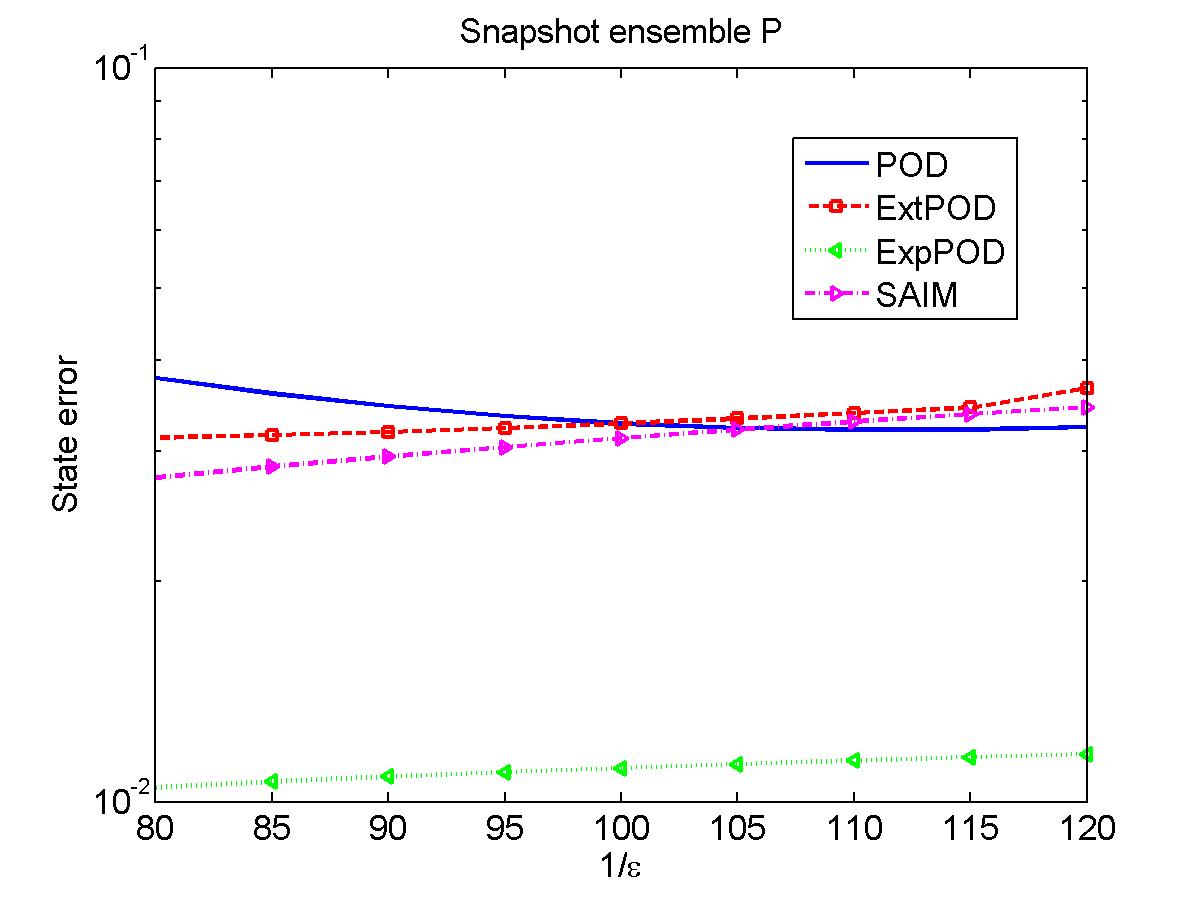}
   \includegraphics[width=0.48\textwidth]{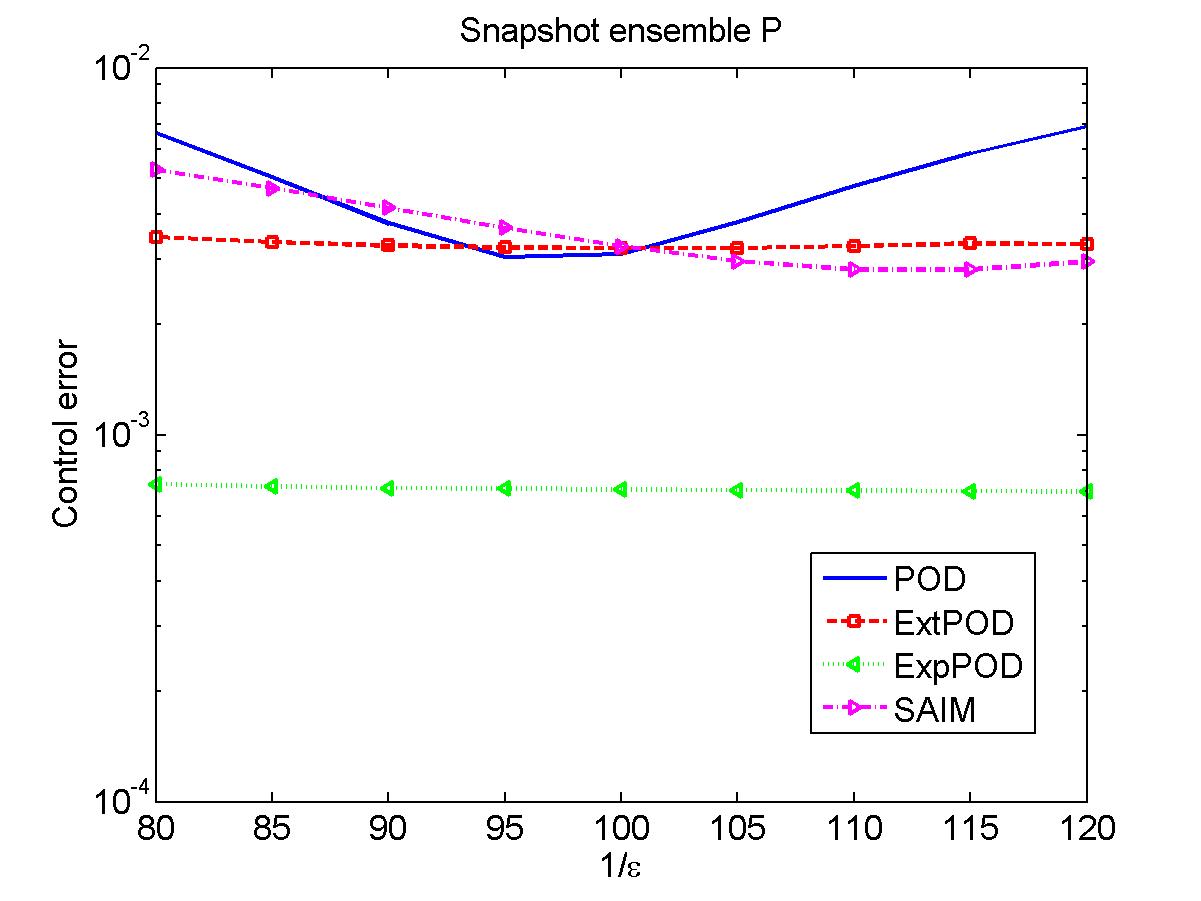}\\
   \includegraphics[width=0.48\textwidth]{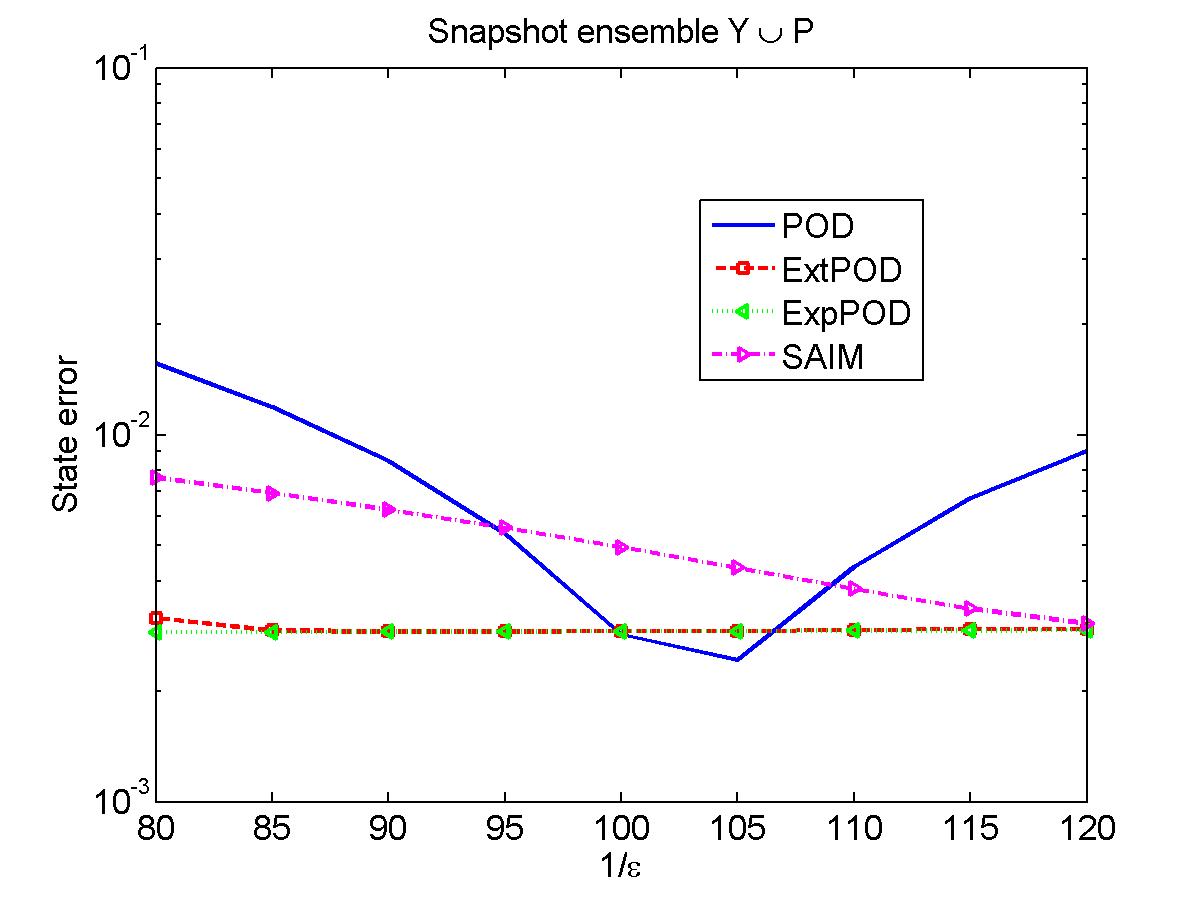}
   \includegraphics[width=0.48\textwidth]{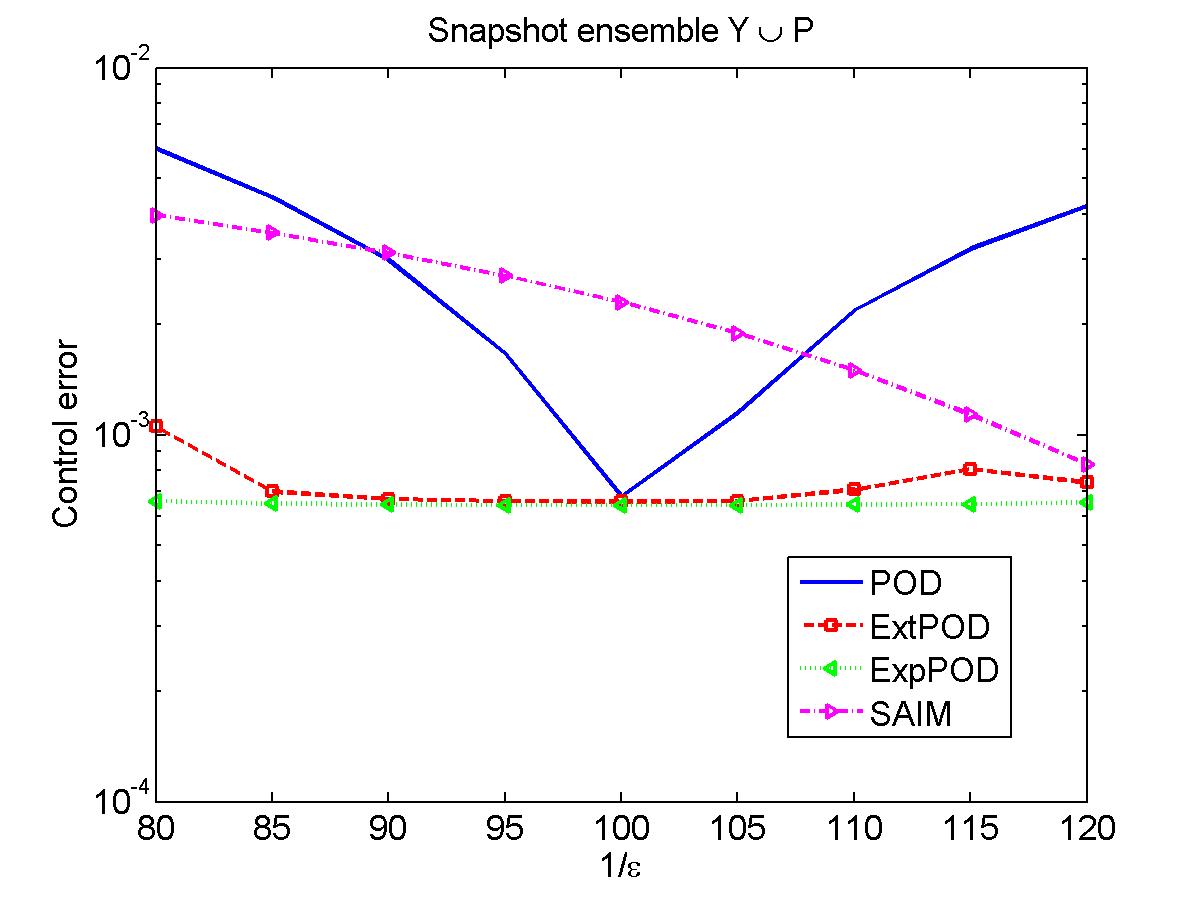}
   \end{array}
      $  
   \caption{The error for the state(\textit{left column}) and the control(\textit{right column}) with 9 POD bases}
   \label{EX2_F:4}
\end{figure}

Now, we compare the low-order solutions obtained by 3 different snapshot sets. We take the numerical solution as the benchmark. In Figure~\ref{EX2_F:4}, we present the error for the state on the left column and the error for the control on the right column with respect to the diffusion term using 9 POD bases. The control approximated with the POD bases generated from the state solution is poor because the characteristics of the control are totally different from the state solution. Similarly, the error in the state approximated by the snapshots of the adjoint is higher. Thus, the choice of the snapshot ensemble affects the approximation depending on whether it contains information about the term which will be approximated or not. The best result for the control is derived with the snapshot set $P$ and $Y \cup P$. The inclusion of the adjoint information in the POD basis generation step improves the performance of the method, because the relation between the adjoint and the control is determined through the optimality condition (\ref{s4_c}). In addition, a good approximation to the control influences the state solution directly due to acting on the right-hand side of (\ref{s2_a}). Therefore, the snapshot ensemble $Y \cup P$ leads to the smallest error. For the state solution, the snapshot set $Y$ and $Y \cup P$ outperforms the snapshot set $P$. It is because the former set contains information about state and the latter offers a good approximation to the control which influences the state solution directly due to the optimality condition. 

We investigate the differences of the POD bases in terms of accuracy depicted in Figure~\ref{EX2_F:4}. We observe that the baseline, ExtPOD and SAIM POD obtained from the snapshot set $Y$ fails to predict the control in case of parameter perturbations since the error remains the same. Instead, ExpPOD improves the approximation. SAIM improves the quality of the perturbed state solution obtain from the snapshot set $Y$ and $Y \cup P$ for larger perturbations. Similar observations are also valid for the control with the snapshot set $P$ and $Y \cup P$. For this example, the best results for both of the state and the control are achieved through the snapshot set $Y \cup P$ with ExtPOD and ExpPOD. As we increase the number of POD bases, ExtPOD becomes superior.

\begin{figure}[H]
   \centering
   $
   \begin{array}{cc} 
   \includegraphics[width=0.48\textwidth]{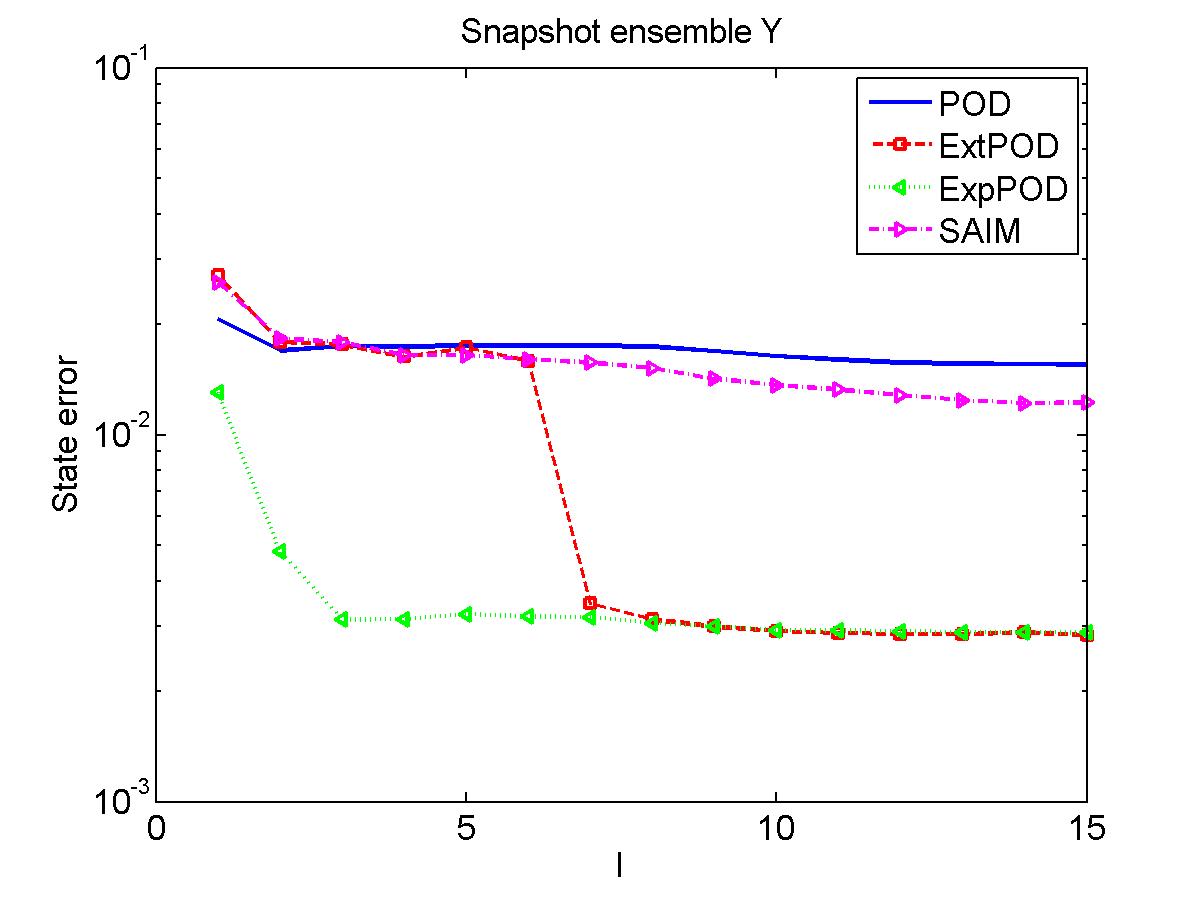}
   \includegraphics[width=0.48\textwidth]{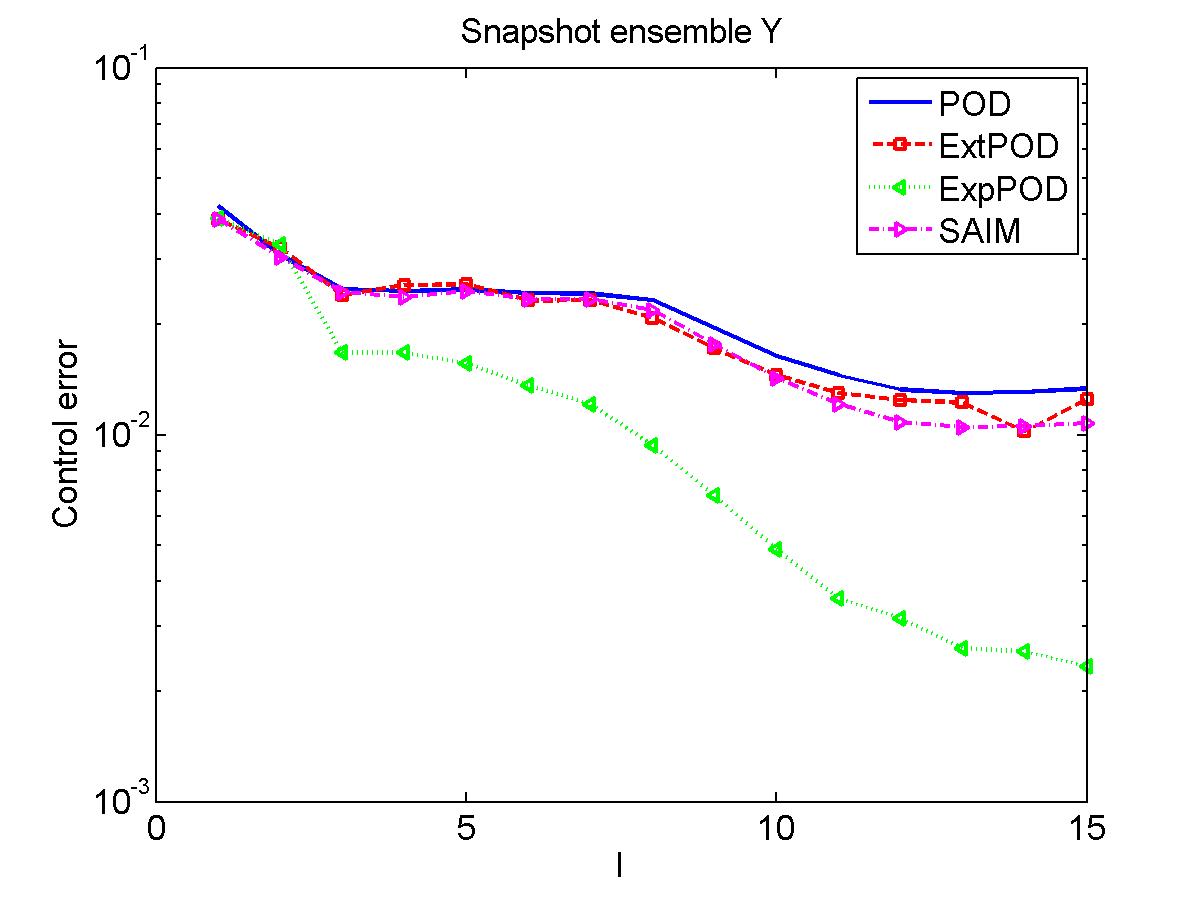}\\
   \includegraphics[width=0.48\textwidth]{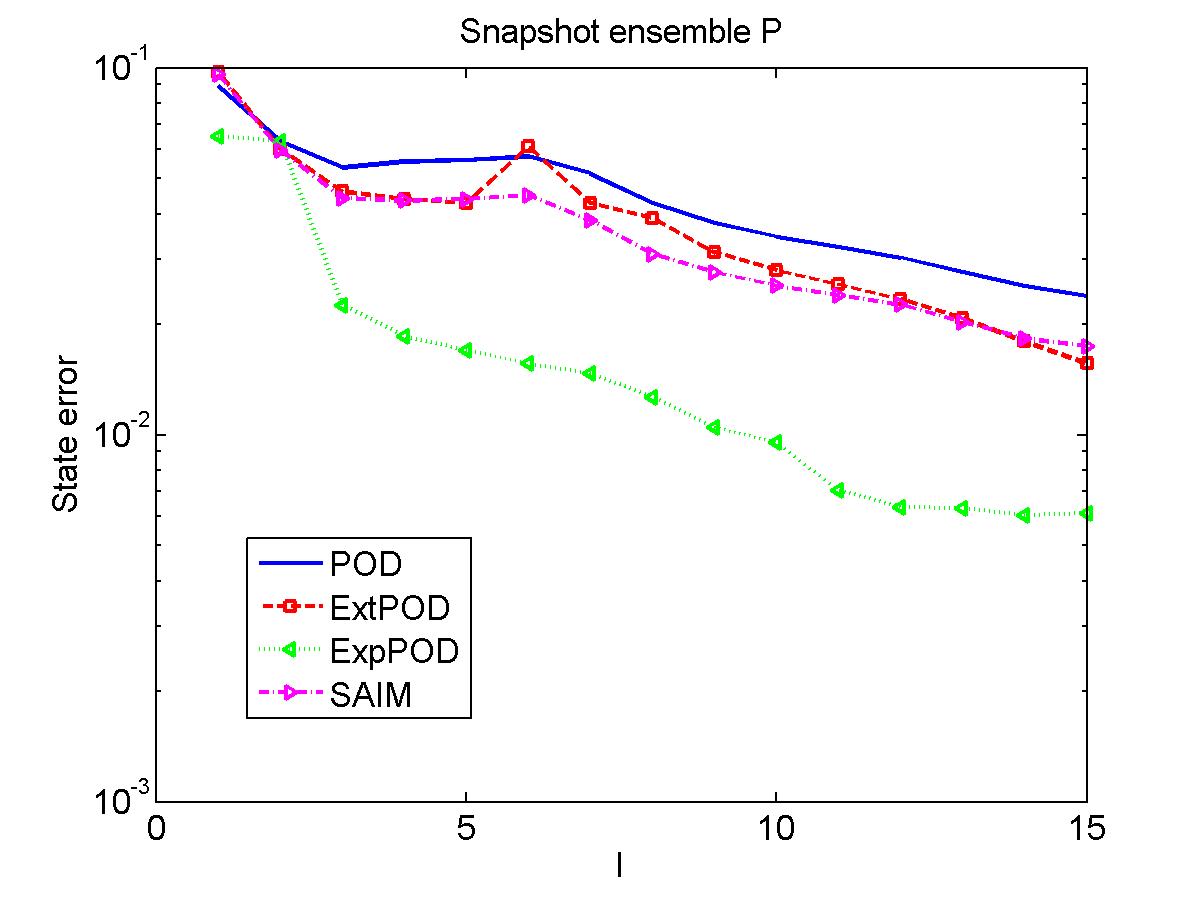}
   \includegraphics[width=0.48\textwidth]{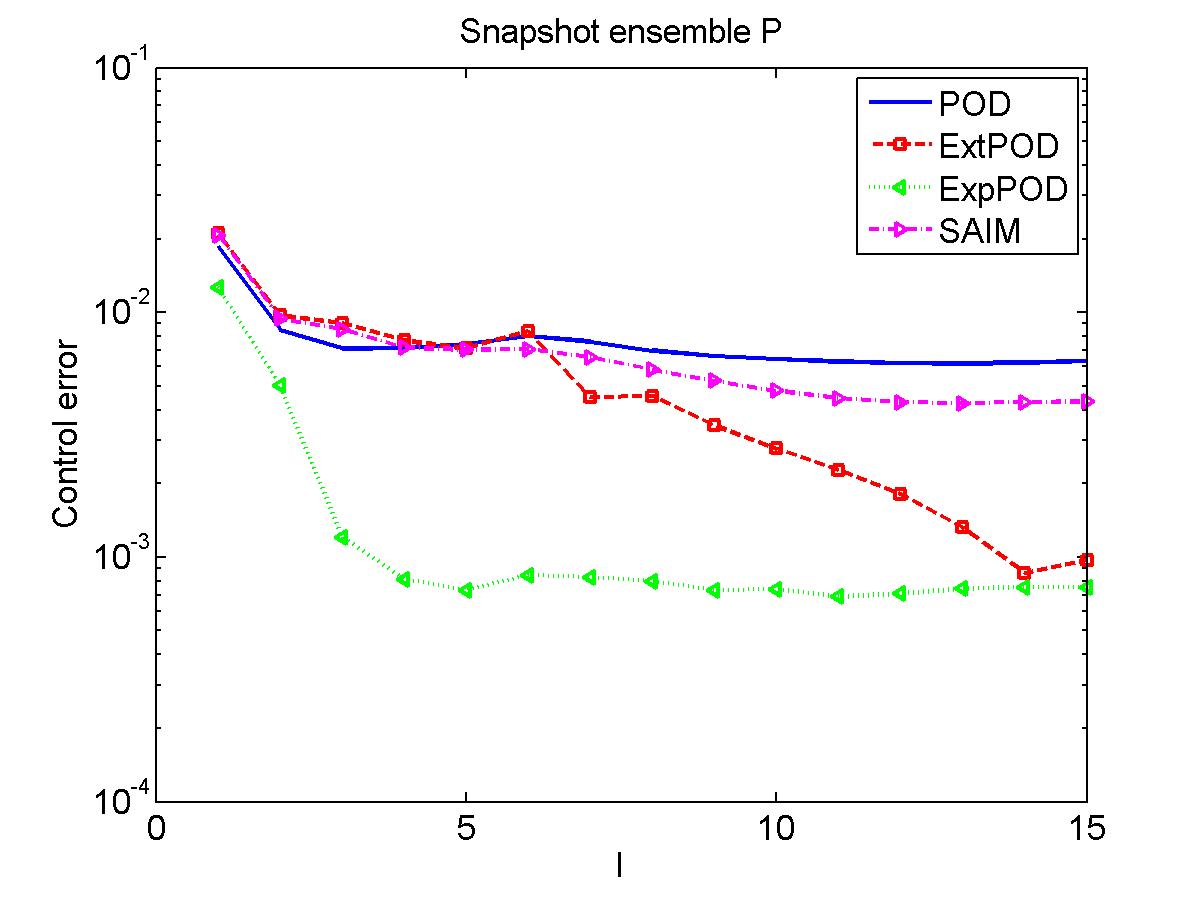}\\
   \includegraphics[width=0.48\textwidth]{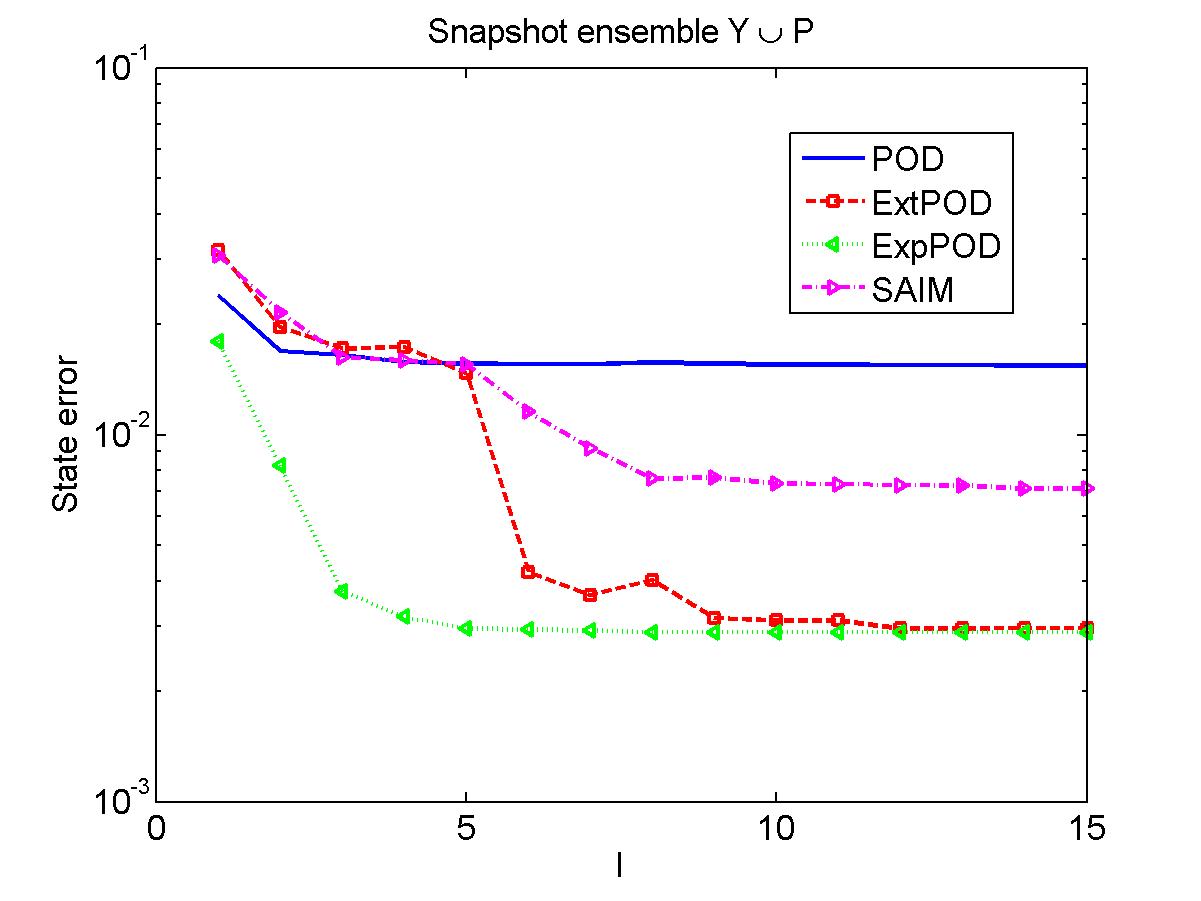}
   \includegraphics[width=0.48\textwidth]{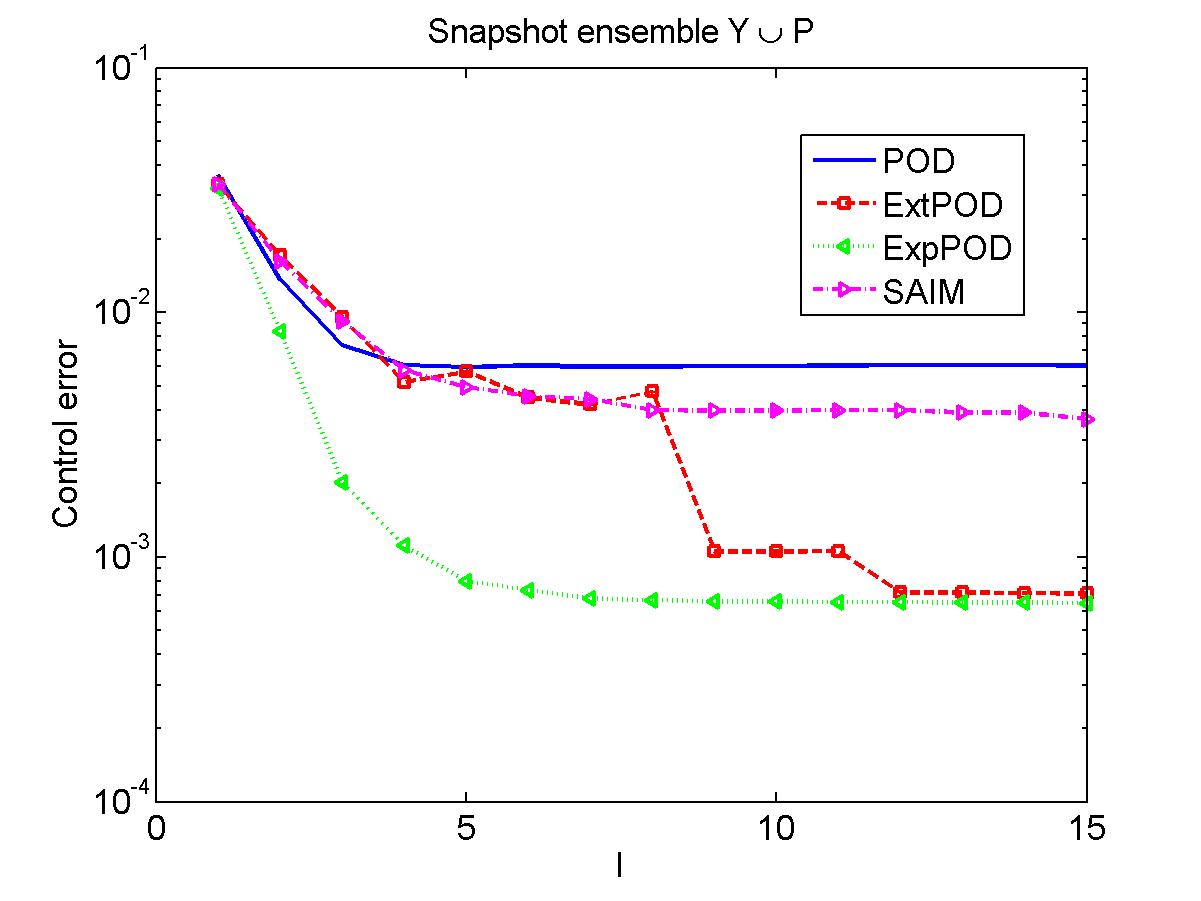}
   \end{array}
   $
   \caption{The error for the state(\textit{left column}) and the control(\textit{right column}) for $\epsilon=1/80$}
   \label{EX2_F:5}
\end{figure}

In Figure~\ref{EX2_F:5}, we present the results for negative changes in $\epsilon$, particularly for $\epsilon=1/80$. All methods surpass the baseline POD. For the state approximation, with the snapshot set $Y$ and $Y \cup P$, ExtPOD approximates the state better than SAIM worse than ExpPOD. However, it catches the ExpPOD as we increase the number of POD bases. On the other hand, for the control, with the snapshot set $P$ and $Y \cup P$, the smallest error is achieved. While ExtPOD offers better results than SAIM, it cannot beat ExpPOD.

\section{Conclusions}\label{conclusion}
We aim to improve the robustness of the POD method by extrapolating or expanding the POD basis using the sensitivity information. In addition, we generate a new basis using SAIM, which requires the snapshots associated to two different parameters $\mu_1$ and $\mu_2$. If one does not have these snapshots, then SAIM would be expensive especially for nonlinear problems. We note that it outperforms the baseline POD for large perturbations. To derive the sensitivities, we observe that CSE method is more time consuming than FD approximation. The snapshot ensemble $Y \cup P$ and ExpPOD leads to the smallest error. The error is almost equal to the one in the ExtPOD as we increase the number of POD bases, so ExtPOD becomes preferable. In a future work, we will be interested in the local improvements to the POD solution of the optimal control of Burgers' equation.




\end{document}